\documentclass[journal]{IEEEtran}
\usepackage{amsmath}
\usepackage{amssymb}
\usepackage[dvips]{graphicx}
\DeclareGraphicsExtensions{.eps}
\usepackage{cite}
\usepackage{multirow}
\usepackage{times}

\usepackage{fancyhdr}
\begin{document}
\title{Factored solution of nonlinear equation systems}
\author{Antonio G\'omez-Exp\'osito, {\it Fellow, IEEE}
\thanks{A. G\'omez-Exp\'osito is with the Department of Electrical Engineering, University of Seville,
Spain; e-mail: age@us.es. 
}}

\maketitle

\begin{abstract}
This article generalizes a recently introduced procedure to solve nonlinear systems of equations, radically departing from the conventional Newton-Raphson scheme. The original nonlinear system is first unfolded into three simpler components: 1) an underdetermined linear system; 2) a one-to-one nonlinear mapping with explicit inverse; 
and 3) an overdetermined linear system. Then, instead of solving such an augmented system at once, a two-step procedure is proposed in which two equation systems are solved at each iteration, one of them involving the same symmetric matrix throughout the solution process. The resulting factored algorithm converges faster than Newton's method and, if carefully implemented, can be computationally competitive for large-scale systems. 
It can also be advantageous for dealing with infeasible cases.
\end{abstract}


%\IEEEpeerreviewmaketitle

%\pagestyle{myheadings}
%\thispagestyle{plain}
%\markboth{ANTONIO G\'OMEZ-EXP\'OSITO}{FACTORED SOLUTION OF NONLINEAR SYSTEMS}

\section{Introduction}

%\PARstart{E}{fficiently}
Efficiently and reliably solving nonlinear equations is of paramount importance in physics, engineering, operational research and many other disciplines in which the need arises to build detailed mathematical models of real-world systems, all of them nonlinear in nature to a certain extent \cite{more,beers,judd,chua}. Moreover, large-scale systems are always sparse, which means that the total number of additive terms in a coupled set of $n$ nonlinear functions in $n$ variables is roughly O($n$).

According to John Rice, who coined the term {\it mathematical software} in 1969, ``solving systems of nonlinear equations is perhaps the most difficult problem in all of numerical computations" \cite{rice}. This surely explains why so many people have devoted so much effort to this problem for so long. 

Since the 17th century, the reference method for the solution of nonlinear systems is Newton-Raphson's (NR) iterative procedure, which locally approximates the nonlinear functions by their first-order Taylor-series expansions. Its terrific success stems from its proven quadratic convergence (when a sufficiently good initial guess is available), moderate computational expense, provided the Jacobian sparsity is fully exploited when dealing with large systems, and broad applicability to most cases in which the nonlinear functions can be analytically expressed \cite{ortega}.

For large-scale systems, by far the most time-consuming step lies in the computation of the Jacobian and its $LU$ factors \cite{tinney}. This has motivated the development of quasi-Newton methods, which make use of approximate Jacobians usually at the expense of more iterations \cite{stott}. Well-known examples in this category are the chord method, where the Jacobian is computed only once, and the secant method \cite{secant}, which approximates the Jacobian through finite differences (no explicit derivatives are required). Broyden's method is a generalization of secant method which carries out rank-one updates to the initial Jacobian \cite{broyden}.

Another category of related methods, denoted as inexact Newton methods, performs approximate computations of Newton steps, frequently in combination with pre-conditioners \cite{kelley}.

More sophisticated higher-order iterative methods also exist, achieving super-quadratic convergence near the solution at the expense of more computational effort. These include Halley's cubic method.

When the initial guess is not close enough to the solution or the functions to be solved exhibit acute non-convexities in the region of interest, the NR method works slowly or may diverge quickly. This has motivated the development of so-called globally convergent algorithms which, for any initial guess, either converge to a root or fail to do so in a small number of ways \cite{ortega,kelley}. Examples of globally convergent algorithms are line search methods, 
continuation/homotopy methods \cite{gina,chao}, such as Levenberg-Marquardt methods, which circumvent Newton's method failure caused by a singular or near singular Jacobian matrix, 
and trust region methods.

Other miscellaneous methods, more or less related to Newton's method, are based on polynomial approximations (e.g. Chebyshev's method), solution of ordinary differential equations  
(e.g. Davidenko's method, a differential form of Newton's method) 
or decomposition schemes (e.g. Adomian's method \cite{adm}).

While each individual in such a plethora of algorithms may be most appropriate for a particular niche application, the standard NR method still remains the best general-purpose candidate trading off simplicity and reliability, particularly when reasonable initial guesses can be made \cite{deuf,kelley2}

A feature shared by most existing methods for solving nonlinear equations is that the structure of the original system is kept intact throughout the solution process. In other words, the algorithms are applied without any kind of preliminary transformation, even though it has long been known that by previously rearranging nonlinear systems better convergence can be achieved. According to \cite{judd} (pp. 174-176), ``the general idea is that a global nonlinear transformation may create an algebraically equivalent system on which Newton's method does better because the new system is more linear. Unfortunately, there is no general way to apply this idea; its application will be problem-specific''.

This article explores such an idea through a new promising perspective, based on unfolding the nonlinear system to be solved by identifying distinct nonlinear terms, each deemed a new variable. This leads to an augmented system composed of two sets of linear equations, which are coupled through a one-to-one nonlinear mapping with diagonal Jacobian. The resulting procedure involves two steps at each iteration: 1) solution of a linear system with symmetric coefficient matrix; 
2) computation of a Newton-like step. 

The basic idea of factoring the solution of complex nonlinear equations into two simpler problems, linear or nonlinear, was originally introduced in the context of hierarchical state estimation \cite{procieee}, where the main goal was to geographically decompose large-scale least-squares problems. Later on, it has been successfully applied to nonlinear equations with a very particular network structure, such as those arising in power systems analysis and operation \cite{eqconst,flf}. 
 
In this work, the factored solution method, initially derived for overdetermined equation systems, is conceptually re-elaborated from scratch, and generalized, so that it can be applied to efficiently solving broader classes of nonlinear systems of equations.

The paper is structured as follows: in the next section the proposed two-stage solution approach is presented. Then, sections \ref{sec3} and \ref{sec4} introduce canonical and extended forms of nonlinear functions to which the proposed method can be applied. Next, section \ref{sec5} discusses how, with minor modifications, the proposed method can reach different solution points from the same initial guess. Section \ref{sec7} is devoted to the analysis of cases which are infeasible in the real domain, where convergence to complex values takes place. Section \ref{sec8} briefly considers the possibility of safely reaching singular points while section \ref{lstc} closes the paper by showing the behaviour of the proposed method when applied to large-scale test cases. 

\section{Factored solution of nonlinear systems}\label{sec2}

Consider a general nonlinear system, written in compact form as follows:
\begin{equation}\label{e1}
h(x)=p
\end{equation}
where $p\in\mathbb{R}^n$ is a specified vector and $x\in\mathbb{R}^n$ is the unknown vector
\footnote{Mathematicians usually write nonlinear equations as $h(x)=0$, the problem being that of finding the roots of the nonlinear function. In engineering and other disciplines dealing with physical systems, however, it is usually preferable to explicitly use a vector of specified parameters, $p$, since the evolution of $x$ (output) as a function of $p$ (input) is of most interest. When no real solutions exist we say that $p$ is infeasible.}

By applying the NR iterative scheme, a solution can be obtained from an initial guess, $x_0$, by successively solving 
\begin{equation}\label{e1i}
H_k \Delta x_k = \Delta p_k  
\end{equation}
where subindex $k$ denotes the iteration counter, $\Delta x_k=x_{k+1}-x_k$, $H_k$ is the Jacobian of $h(\cdot)$, computed at the current point $x_k$, and $\Delta p_k=p-h(x_k)$ is the mismatch or residual vector.

In essence, the new method assumes that suitable auxiliary vectors, $y$ and $u\in\mathbb{R}^m$, with $m>n$, can be introduced so that the original system (\ref{e1}) can be unfolded into the following simpler problems:
\begin{eqnarray}
Ey &=& p \label{e2}\\
u &=& f(y) \label{e3}\\
Cx &=& u \label{e4}
\end{eqnarray}
where $E$ and $C$ are rectangular matrices of sizes $n\times m$ and $m\times n$, respectively, and vector $f(\cdot)$ comprises a set of one-to-one nonlinear functions, also known as a diagonal nonlinear mapping \cite{sandberg},
\begin{equation}\label{e5}
u_i=f_i(y_i) \quad ; \quad i=1,\dots m 
\end{equation}
each with closed-form inverse,
\begin{equation}\label{e6}
y_i=f^{-1}_i(u_i) %\quad ; \quad i=1,\dots m
\Rightarrow  y=f^{-1}(u)
\end{equation}

By eliminating vector $u$ the above augmented system can also be written in more compact form:
\begin{eqnarray}
Ey &=& p \label{e7}\\
Cx &=& f(y) \label{e8}
\end{eqnarray}
Notice that (\ref{e7}) is a linear underdetermined system whereas (\ref{e8}) is an overdetermined one. Among the infinite solutions to (\ref{e7}) only those exactly satisfying (\ref{e8}) constitute solutions to the original nonlinear system (\ref{e1}). As explained in section \ref{sec3}, many systems can be found in practice where such a factorization is possible, the aim being to reduce $m$ as much as possible. In the limit, if a vector $y$ of size $m=n$ can be found, then solutions to the original nonlinear system will be obtained directly (i.e., without iterating) by sequentially solving the unfolded system of equations (\ref{e2})-(\ref{e4}). But this case will arise only when the set of $n$ equations comprises just $n$ nonlinear distinct terms. Therefore, the need to iterate in the factored method stems from an excess of nonlinear terms ($m>n$) rather than from the nonlinear nature of the original problem.

Obviously, when the remaining auxiliary vector $y$ is eliminated the original `folded' system is obtained in factored form:
\begin{equation}\label{e9}
h(x)=Ef^{-1}(Cx)=p
\end{equation}
This leads to an equivalent expression for the Newton step (\ref{e1i}),
\begin{equation}\label{e9i}
[\underbrace{EF_k^{-1}C}_{H_k}] \Delta x_k = p-Ef^{-1}(Cx_k) = \Delta p_k 
\end{equation}
where $F$ is the trivial Jacobian of $f(\cdot)$. Whether (\ref{e9i}) offers any computational advantage over (\ref{e1i}) will mainly depend on the complexity involved in the computation of $h(\cdot)$ and its Jacobian $H$, compared to their much simpler counterparts $f(\cdot)$ and $F$ (plus the required matrix products), but the convergence pattern will be exactly the same in both cases.

Yet the augmented system (\ref{e7})-(\ref{e8}), arising when the factored form is considered, opens the door to alternative solution schemes which may outperform the `vanilla' NR approach. The scheme developed in the sequel begins by considering the solution of the factored model as a linearly-constrained Least Squares (LS) problem. Then, an auxiliary least-distance problem is formulated aimed at providing improved linearization points at each iteration.

\subsection{Solution of the factored model as an optimization problem}\label{sec2.1}

Finding a solution to the unfolded problem (\ref{e7})-(\ref{e8}) can be shown to be equivalent to solving the following equality-constrained LS problem,
\begin{eqnarray}
{\rm Minimize \ \ } && [y-f^{-1}(Cx)]^T[y-f^{-1}(Cx)]  \label{eqwls} \\
{\rm s.t. \ \ }  &&  p-Ey =0   \nonumber 
\end{eqnarray}
which reduces to minimizing the associated Lagrangian function,
\begin{equation}\label{lag1}
{\cal L}= \frac{1}{2}[y-f^{-1}(Cx)]^T[y-f^{-1}(Cx)] - \mu^T(p - Ey)
\end{equation}

The first-order optimality conditions (FOOC) give rise to the following system:
\begin{eqnarray}\label{FOOC}
y-f^{-1}(Cx) + E^T\mu &=& 0 \nonumber \\
-C^TF^{-T}[y-f^{-1}(Cx) ] &=& 0 \\
p - Ey &=& 0 \nonumber
\end{eqnarray} 
 
Given an estimate of the solution point, $x_k$, we can choose $y_k$ in such a way that (\ref{e8}) is satisfied, i.e.,
\begin{equation}
y_k=f^{-1}(Cx_k) \label{eini}
\end{equation}
Then, linearizing $f^{-1}(\cdot)$ around $x_k$
\[ y-f^{-1}(Cx) \cong  \Delta y_k - F_k^{-1}C\Delta x_k  \]
allows (\ref{FOOC}) to be written in incremental form,
\begin{equation}\label{FOOCinc}
\left[\begin{array}{ccc} I & -F_k^{-1}C & E^T \\ 
-C^TF_k^{-T} & C^TD_kC & 0 \\
E & 0 & 0
\end{array}\right]
\left[\begin{array}{c} \Delta{y}_k \\  \Delta{x}_k \\ \mu
\end{array}\right] =
\left[\begin{array}{c} 0 \\ 0 \\ \Delta p_k
\end{array}\right] 
\end{equation}
where $D_k= F_k^{-T}F_k^{-1}$ is a diagonal matrix with positive elements.

From the above symmetric system, $\Delta{y}_k$ can be easily eliminated, yielding:
\begin{equation}\label{FOOCnoy}
\left[\begin{array}{cc} 0 & C^TF_k^{-T}E^T \\ 
EF_k^{-1}C & -EE^T
\end{array}\right]
\left[\begin{array}{c} \Delta{x}_k \\ \mu
\end{array}\right] =
\left[\begin{array}{c} 0 \\ \Delta p_k
\end{array}\right] 
\end{equation}
or, in more compact form,
\begin{equation}\label{FOOCnoyc}
\left[\begin{array}{cc} 0 & H_k^T \\ 
H_k & -EE^T
\end{array}\right]
\left[\begin{array}{c} \Delta{x}_k \\ \mu
\end{array}\right] =
\left[\begin{array}{c} 0 \\ \Delta p_k
\end{array}\right] 
\end{equation}
If the Jacobian $H_k$ remains nonsingular at the solution point, then $\mu=0$ and the above system reduces to (\ref{e9i}), as expected. Therefore, the Lagrangian-based augmented formulation has the same convergence pattern and converges to the same point as the conventional NR approach.

\subsection{Auxiliary least-distance problem}\label{sec2.2}

The driving idea behind the proposed procedure is to linearize (\ref{FOOC}) around a point which, being closer to the solution than $y_k$, can be obtained in a straightforward manner. 
The resulting scheme will be advantageous over the NR approach if the extra cost of computing the improved linearization point is offset by the convergence enhancement obtained, if any. 

For this purpose, given the latest solution estimate, $x_k$, we consider a simpler auxiliary optimization problem, as follows:
\begin{eqnarray}
{\rm Minimize \ \ } && (y-y_k)^T(y-y_k)  \label{seqwls} \\
{\rm s.t. \ \ }  &&  p - Ey =0   \nonumber 
\end{eqnarray}
with $y_k=f^{-1}(Cx_k)$. The associated Lagrangian function is,
\begin{equation}\label{lags}
{\cal L}_a= \frac{1}{2}(y-y_k)^T(y-y_k) - \lambda^T(p - Ey)
\end{equation}
which leads to the following linear FOOCs:
\begin{equation}\label{FOOCaux}
\left[\begin{array}{cc} I & E^T \\ 
E & 0
\end{array}\right]
\left[\begin{array}{c} \Delta\tilde{y}_k \\ \lambda
\end{array}\right] =
\left[\begin{array}{c} 0 \\ \Delta p_k
\end{array}\right] 
\end{equation}
where $\Delta\tilde{y}_k=\tilde{y}_k-y_k$. Besides being linear, the unknown $x$ is missing in this simpler problem, which can be solved in two phases. First, $\lambda$ is computed by solving:
\begin{equation}
EE^T \lambda= - \Delta p_k \label{lam}
\end{equation}
Then, $\tilde{y}_k$ is simply obtained from 
\begin{equation}
\tilde{y}_k= y_k - E^T\lambda \label{imy}
\end{equation}
By definition, $\tilde{y}_k$ is as close as possible to $y_k$ while satisfying (\ref{e7}).

Next, the FOOCs of the original problem (\ref{FOOC}) are linearized around $\tilde{y}_k$, hopefully closer to the solution than $y_k$,
\[ y-f^{-1}(Cx) \cong  \Delta\tilde{y}_k - \tilde{F}^{-1}[Cx - f(\tilde{y}_k)] \]
leading to,
\begin{equation}\label{FOOCauxinc}
\left[\begin{array}{ccc} I & -\tilde{F}^{-1}C & E^T \\ 
-C^T\tilde{F}^{-T} & C^T\tilde{D}C & 0 \\
E & 0 & 0
\end{array}\right]
\left[\begin{array}{c} \Delta\tilde{y}_k \\  {x}_{k+1} \\ \mu
\end{array}\right] =
\left[\begin{array}{c} - \tilde{F}^{-1}f(\tilde{y}_k) \\ C^T\tilde{D}f(\tilde{y}_k) \\ 0
\end{array}\right] 
\end{equation}

Once again, eliminating $\Delta\tilde{y}_k$ yields,
\begin{equation}\label{FOOCauxnoy}
\left[\begin{array}{cc} 0 & C^T\tilde{F}^{-T}E^T \\ 
E\tilde{F}^{-1}C & -EE^T
\end{array}\right]
\left[\begin{array}{c} {x}_{k+1} \\ \mu
\end{array}\right] =
\left[\begin{array}{c} 0 \\ E\tilde{F}^{-1}f(\tilde{y}_k)
\end{array}\right] 
\end{equation}
or, in more compact form,
\begin{equation}\label{FOOCauxnoyc}
\left[\begin{array}{cc} 0 & \tilde{H}^T \\ 
\tilde{H} & -EE^T
\end{array}\right]
\left[\begin{array}{c} {x}_{k+1} \\ \mu
\end{array}\right] =
\left[\begin{array}{c} 0 \\ E\tilde{F}^{-1}f(\tilde{y}_k)
\end{array}\right] 
\end{equation}
The above linear system, to be compared with (\ref{FOOCnoyc}), provides the next value ${x}_{k+1}$, which replaces ${x}_{k}$ in the next iteration. Moreover, as happened with the solution of (\ref{FOOCnoyc}), so long as the Jacobian $\tilde{H}$ remains nonsingular, $\mu=0$ and ${x}_{k+1}$ can be obtained with less computational effort from:
\begin{equation}
\tilde{H}{x}_{k+1} = E\tilde{F}^{-1}f(\tilde{y}_k) \label{newx0}
\end{equation}

\subsection{Two-step solution procedure}\label{sec2.3}

The two-step algorithm, arising from the proposed factored representation of nonlinear systems, can then be formally stated as follows: 

\vspace{10pt}
\framebox{\begin{minipage}{0.45\textwidth}

\noindent{\bf Step 0: Initialization.} 
Initialize the iteration counter ($k=0$). Let $x_k$ be an initial guess and $y_k=f^{-1}(Cx_k)$.

\noindent{\bf Step 1: Auxiliary least-distance problem.} 
First, obtain vector $\lambda$ by solving the system
\begin{equation} (EE^T)\lambda=p - E{y}_{k} \label{newl}\end{equation}
and then compute $\tilde{y}_k$ from,
\begin{equation}
\tilde{y}_k= y_k + E^T\lambda \label{newy}
\end{equation}

\noindent{\bf Step 2: Non-incremental computation of $x$.}   
Solve for ${x}_{k+1}$ the system,
\begin{equation}
\tilde{H}{x}_{k+1} = E\tilde{F}^{-1}f(\tilde{y}_k) \label{newx}
\end{equation}
where $\tilde{H} = E\tilde{F}^{-1}C$ is the factored Jacobian computed at $\tilde{u}_k=f(\tilde{y}_k)$.
Then update ${y}_{k+1}=f^{-1}(C {x}_{k+1})$. If $||x_{k+1}-x_k||$ (or, alternatively $||p - E{{y}_{k+1}}||$) is small enough, then stop. Otherwise set $k=k+1$ and go back to Step 1.

\end{minipage}}\vspace{10pt}

As explained above, step 1 constitutes a linearly-constrained least-distance problem, yielding a vector $\tilde{y}_k$ which, being as close as possible to $y_k$, satisfies (\ref{e7}). As a feasible solution is approached, $\lambda\rightarrow 0$ and $\tilde{y}_k\rightarrow y_k$. Notice that the sparse symmetric matrix $EE^T$ needs to be factorized only once, by means of the efficient Cholesky algorithm (in fact, only its upper triangle is needed). Therefore, the computational effort involved in step 1 is very low if Cholesky triangular factors are saved during the first iteration. Moreover, the vector $p - E{y}_{k}$ is available from the previous iteration if it is computed to check for convergence. 

It is worth noting that step 2 is expressed in non-incremental form. It directly provides improved values for $x$ with less computational effort than that of a conventional Newton step (\ref{e9i}), which may well offset the moderate extra burden of step 1. 

In \cite{eqconst,flf} the two-step procedure based on the factored model was originally developed from a different perspective, related with the solution of overdetermined systems in state estimation (i.e. filtering out noise from a redundant set of measurements). 

\subsection{Factored versus conventional Newton's method} \label{sec2.4}

The factored scheme can be more easily compared with NR method if step 2 is written in incremental form. For this purpose, the term $\tilde{H}{x}_{k}$ is subtracted from both terms of (\ref{newx}). Considering that $f(y_k)=Cx_k$, this leads to,
\begin{equation}
\tilde{H}({x}_{k+1}-x_k) = E\tilde{F}^{-1}[f(\tilde{y}_k) - f(y_k)] \label{newxI}
\end{equation} 
Notice that the above expression, being equivalent to (\ref{newx}), is less convenient from the computational point of view owing to the extra operations involved.

Next, the Taylor's expansion around $\tilde{y}_k$ of the set of nonlinear functions $f(y)$ is rearranged as follows, for  $y=y_k$:
\begin{equation}
f(\tilde{y}_k)-f(y_k)= \tilde{F}(\tilde{y}_k-y_k) - \mathcal{R}(y_k,\tilde{y}_k)  
\label{f_taylor}
\end{equation}
where the remainder $\mathcal{R}(y_k,\tilde{y}_k)$ comprises the second- and higher-order terms of the polynomial expansion,
\begin{equation}
\mathcal{R}(y_k,\tilde{y}_k) 
= \sum_{j=2}^{\infty}\frac{\tilde{F}_j}{j!}[{\rm diag\ }(y_k-\tilde{y}_k)]^j e
\label{rem}
\end{equation}
$\tilde{F}_j$ is the diagonal matrix containing the $j$-th derivatives of  $f(y)$, computed at $y=\tilde{y}_k$, and $e$ is a vector of 1's. 

Replacing (\ref{f_taylor}) into (\ref{newxI}), and taking into account that $E\tilde{y}_k=p$, yields, 
\begin{equation}
\tilde{H}\Delta x_k = \Delta p_k - E\tilde{F}^{-1}\mathcal{R}(y_k,\tilde{y}_k)  
\label{newxTI}
\end{equation}

Comparing the resulting incremental model (\ref{newxTI}) with that of the conventional Newton's method (\ref{e9i}), the following remarks can be made:
\begin{enumerate} 
\item If the current point, $x_k$, is so close to the solution that $||\tilde{y}_k-y_k||<<1$, then $\mathcal{R}$ will be negligible compared to $\Delta p_k$. In these cases, the only difference of the factored method with respect to NR lies in the use of the updated Jacobian, $\tilde{H}$, rather than $H_k$, and the local convergence rate of the factored algorithm will be of the same order as that of Newton's method. 
\item If the current point, $x_k$, is still so far away from the solution that $||\tilde{y}_k-y_k||>1$, then $\mathcal{R}$ will dominate the right-hand side in (\ref{newxTI}). In these cases, the behaviour of the factored scheme can be totally different from that of Newton's method, which fully ignores $\mathcal{R}$ when computing $x_{k+1}$. 
\item As expected, if step 1 is omitted (i.e., $\tilde{y}_k=y_k$), the factored scheme reduces to the standard Newton's method. 
\end{enumerate}

Let $x_\infty$ and $f(y_\infty)=Cx_\infty$ denote the solution point. Further insight into Remark (i) above, regarding local convergence rate, can be gained by subtracting $\tilde{H}x_\infty$ from both sides of (\ref{newx}) and rearranging, which leads for the factored method to:
\begin{equation}
\tilde{H}(x_\infty - x_{k+1}) = E\tilde{F}^{-1}\mathcal{R}(y_\infty,\tilde{y}_k)  
\label{quadcF}
\end{equation}
Manipulating (\ref{e9i}) in a similar fashion yields, for the NR method:
\begin{equation}
H_k(x_\infty - x_{k+1}) = EF_k^{-1}\mathcal{R}(y_\infty,y_k)  
\label{quadcNR}
\end{equation}
The last two expressions confirm that both the NR and the factored schemes converge quadratically near the solution (third- and higher-order terms are negligible in $\mathcal{R}$). Moreover, if $\tilde{y}_k$ is closer than $y_k$ to the solution $y_\infty$, then the factored method converges faster than Newton's.

In summary, from the above discussion it can be concluded that the local behaviour of the proposed factored algorithm is comparable to Newton's method (quadratic convergence rate), while the global behaviour (shape of basins of attractions) will be in general quite different. The fact that step 1 tries to minimize the distance between $\tilde{y}_k$ and $y_k$ surely explains the much wider basins of attractions found in practice for the factored method, but there is no way to predict beforehand the global behaviour of any iterative algorithm in the general case (the reader is referred to Section 3, ``Aspects of Convergence Analysis'', in \cite{rheinboldt}, for a succinct but clarifying discussion on this issue).

\section{Application to canonical forms}\label{sec3}

The factored model and, consequently, the two-stage procedure presented above can be applied in a straightforward manner to a wide range of nonlinear equation systems, which are or can be directly derived from the following canonical forms.
Small tutorial examples will be used to illustrate the ideas. None of those examples are intended to assess the potential computational benefits of the factored solution process, but rather to show the improved convergence pattern over the NR method (the reader is referred to Section \ref{lstc}, where very large equation systems arising in the steady-state solution of electrical power systems are solved by both methods and solution times are compared).

\subsection{Sums of single-variable nonlinear elementary functions}\label{sec3.1}
The simplest canonical form of nonlinear systems arises when they are composed of sums of nonlinear terms, each being an elementary invertible function of just a single variable. Mathematically, each component of $h(x)$ should have the form,
\[  h_i(x) = \sum_{j=1}^{m_i} c_{ij} h_{ij}(x_k)  \]
where $c_{ij}$ is any real number and $h^{-1}_{ij}(\cdot)$ can be analytically obtained.
In this case, for each different nonlinear term, $h_{ij}(x_k)$, a new variable
\[  y_{ij}= h_{ij}(x_k)  \]
is added to the auxiliary vector $y$, keeping in mind that repeated definitions should be avoided. This leads to the desired underdetermined linear system to be solved at the first step: 
\[  h_i(x) = \sum_{j=1}^{m_i} c_{ij} y_{ij}  \]
while matrix $C$ of the overdetermined linear system $u=Cx$ is trivially obtained from the definition of $u$:
\[  u_{ij}=h^{-1}_{ij}(y_{ij}) = x_k \]

\vspace{3pt}
\noindent{\bf Example 1:}

\begin{figure}[!h]
\begin{center}
\noindent\includegraphics[width=0.45\textwidth]{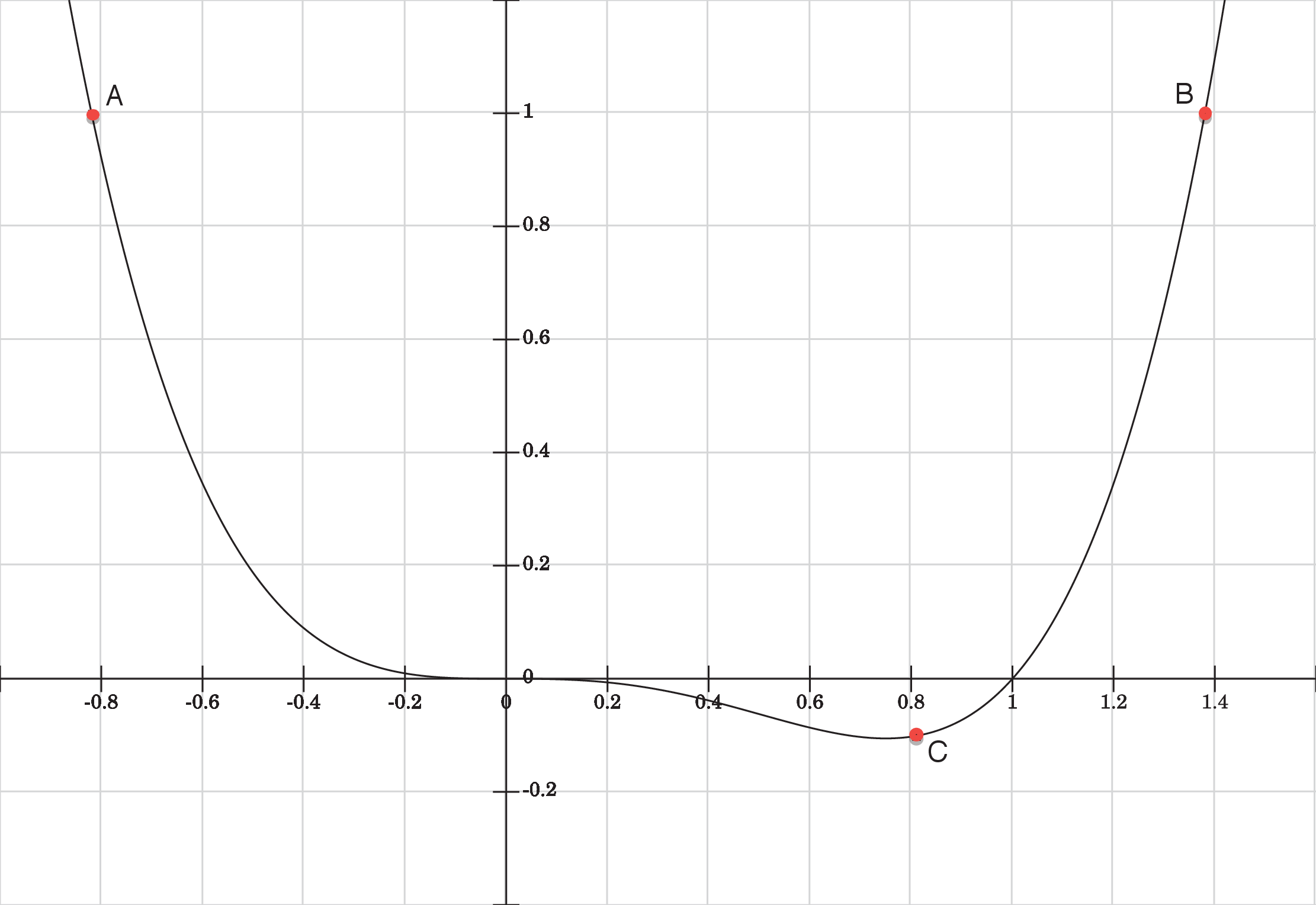}
\caption{Nonlinear function $x^4-x^3$ used in {\bf Example 1}. It crosses $p=1$ at points A and B ($x=-0.8192$ and $x=1.3803$, respectively).} \label{fex1}
\end{center}
\end{figure}

Consider the simple nonlinear function,
\[  p=x^4-x^3 \]
represented in Fig. \ref{fex1}.  Introducing two auxiliary variables, 
\[  y_1=x^4 \quad ; \quad y_2=x^3 \]
leads to,
\[  p=y_1-y_2 \]

Therefore, with the above notation, the involved matrices are:
\[  
E=\left(  \begin{array}{cc} 1 & -1 \end{array} \right) \,\, ; \,\,
C=\left( \begin{array}{c} 1  \\ 1 \end{array}\right) \] 
\[ u=f(y)= \left(  \begin{array}{c} \sqrt[4]{y_1} \\  \sqrt[3]{y_2} \end{array} \right) 
\,\, ; \,\,  F^{-1}= \left(\begin{array}{cc} 4u_1^3 & 0 \\ 0  & 3u_2^2
\end{array}\right) 
\]

For $p=1$ there are two solutions (points A and B in figure \ref{fex1}). Table \ref{t1ex1} shows, 
for different starting points, the number of iterations required to converge (convergence tolerance in all examples: $||\Delta x||_1< 0.00001$) by both the NR and the proposed factored procedure. Notice 
that the factored procedure converges much faster than the NR method. The farther $x_0$ from the solution point, the better the behaviour of the proposed procedure compared to Newton's method.  

\begin{table}[hbt]
\caption{Number of iterations required by both NR and factored procedures to converge from arbitrary starting points for the nonlinear function of Example 1, with $p=1$. In parenthesis, A or B indicates which point is reached in each case.}\label{t1ex1}
  \centering
  \begin{tabular}{|r||r|r|}
    \hline
    % after \\: \hline or \cline{col1-col2} \cline{col3-col4} ...
    $x_0$ & NR & Factored \\ \hline
    30  & 16 (B) & 6 (B)\\
    10  &  12 (B)& 6 (B)\\  % 5
     5 & 9 (B)&  5 (B)\\     % 3
     1 & 7 (B)&  4 (B)\\
    0.9  & 9 (B)& 5 (B)\\ %
    0.8  & 13 (B) & 5 (B)\\
    0.5  & 10 (A)& 6 (B)\\
    0  &  Fails &  6 (B)\\
  $-0.5$  &  6 (A)  &  7 (B)\\
    \hline
  \end{tabular}
\end{table}

For values of $x_0<0.75$ (the minimum of the function, where the slope changes its sign), the NR method converges to point A. On the other hand, the factored scheme always converges to point B no matter which initial guess is chosen. This issue is discussed in more detail in section \ref{sec5}, where it is also explained how to reach other solution points.  As expected, the NR method fails for $x_0=0$ (null initial slope), whereas the factored solution does not suffer from this limitation (the updated Jacobian $\tilde{H}$ corresponding to the first value $\tilde{y}_k$ is not singular). 

Finally, the components of vector $\lambda$ (and obviously $\mu$) are always null at the solution point, indicating that a feasible solution has been found in all cases (see section \ref{sec7}).

\vspace{3pt}
\noindent{\bf Example 2:}

In this example, the following periodic function will be considered,
\[  p=\sin{x} + \cos{x} \]
for which two auxiliary variables are needed, 
\[  y_1=\sin{u_1} \quad ; \quad y_2=\cos{u_2} \]
with $u_1=u_2=x$.
The relevant matrices arising in this case are:
\[  
E=\left(  \begin{array}{cc} 1 & 1 \end{array} \right) \,\, ; \,\,
C=\left( \begin{array}{c} 1  \\ 1 \end{array}\right) \] 
\[ u=f(y)= \left(  \begin{array}{c} \arcsin{y_1} \\  \arccos{y_2} \end{array} \right) 
\,\, ; \,\,  F^{-1}= \left(\begin{array}{cc} \cos{u_1} & 0 \\ 0  & -\sin{u_2}
\end{array}\right) 
\]

\begin{figure}[hbt]
\begin{center}
\noindent\includegraphics[width=0.45\textwidth]{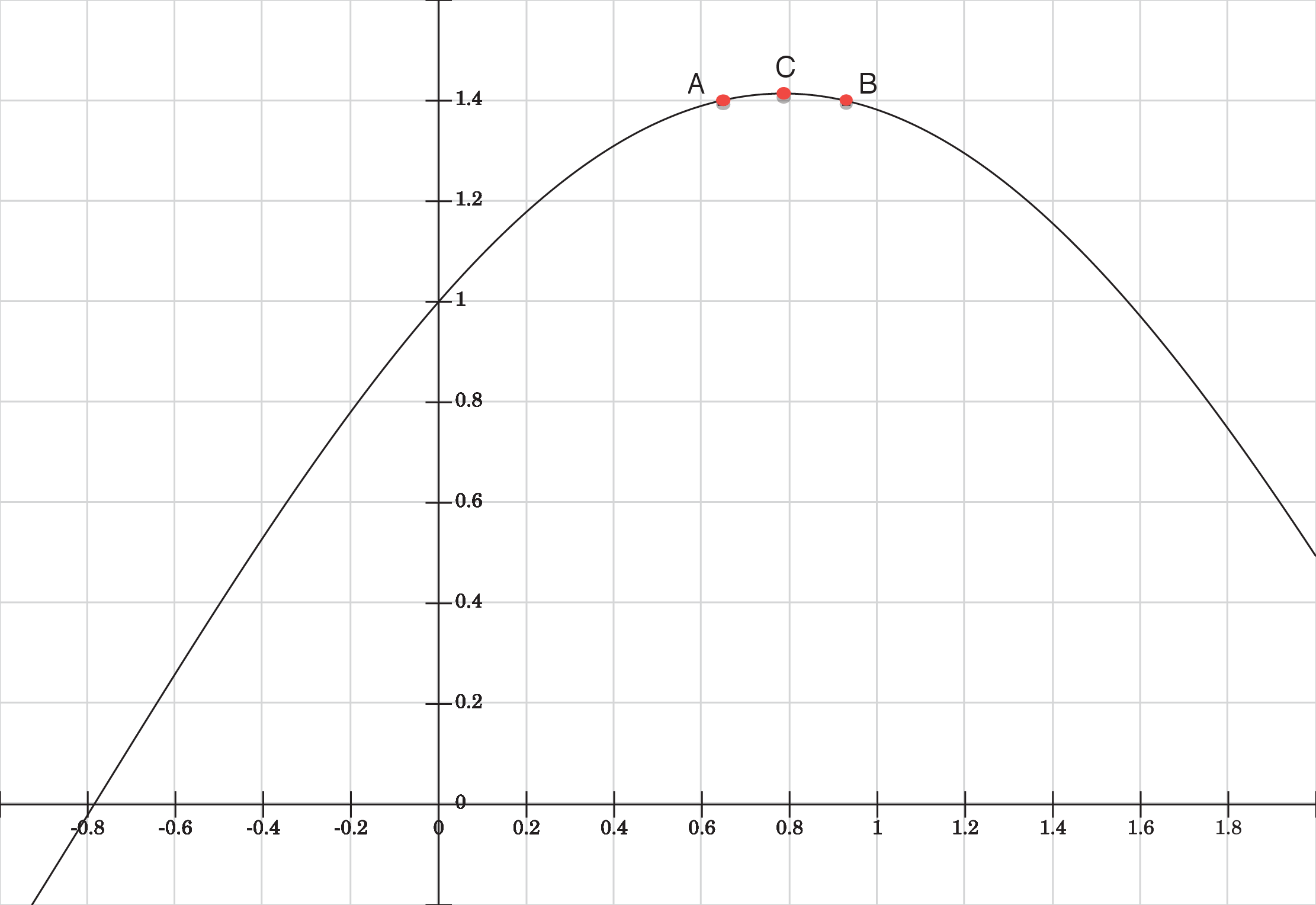}
\caption{Nonlinear function for Example 2} \label{fex2}
\end{center}
\end{figure}

For $p=1.4$, the two solutions closest to the origin are $x=0.6435$ and $x=0.9273$ (points A and B in Fig. \ref{fex2}).

Table \ref{t1ex2} shows, for different starting points, the number of iterations required to converge by both methods and the final solution reached. 
In this particular example ($p=1.4$), the NR method sometimes outperforms the factored scheme in terms of number of iterations. However, note that, for starting points far away from the origin, the NR method may converge to `remote' periodic solutions (shown in boldface), whereas the factored scheme always converges to points A or B (this issue if further discussed in section \ref{sec5}). Other advantages of the factored method will become apparent for infeasible values of $p$ (see the same example, with $p=1.5$, in section \ref{sec7}).

\begin{table}[hbt]
  \caption{Example 2: No. of iterations and solution points for  $p=1.4$}\label{t1ex2}
  \centering
  \begin{tabular}{|r||c|r||c|r|}
    \hline
    $x_0$ & \multicolumn{2}{c||}{NR} & \multicolumn{2}{c|}{Factored} \\ \cline{2-5}
     &  Iter & $x$ &  Iter & $x$ \\ \hline
     10     &  5  & 0.6435   &   4   & 0.9273  \\
     5     &  6  & {\bf 6.9267}  &   8   & 0.6435  \\
     1     &  4  & 0.9273   &   4  &  0.9273  \\
     0     &  6  &  0.6435  &   7  &   0.6435  \\
    $-1$     &  7  &  0.6435  &  8  &   0.6435  \\
     $-5$     &  6  & $-${\bf 55.6214}   &   7   & 0.9273  \\
   $-10$     &  7  & $-${\bf 11.6391}   &   8   & 0.9273  \\
    \hline
  \end{tabular}
\end{table}

In all cases, the final solution provided by the factored scheme is real (or, more precisely, its imaginary component is well below the convergence tolerance). However, complex intermediate solutions are obtained in some cases, owing to the fact that $\arcsin{y}$ and $\arccos{y}$ return complex values for $|y|>1$ (this is discussed below in section \ref{sec7}). 

\subsection{Sums of products of single-variable power functions}\label{sec3.2}

Another canonical form to which the factored method can be easily applied arises when the nonlinear equations to be solved are the sum of products of nonlinear terms, each being an arbitrary power of a single variable. Mathematically, each component of $h(x)$ has the form,
\[  h_i(x) = \sum_{j=1}^{m_i} c_{ij} \prod_{k=1}^{n_j} x_k^{q_k}  \]
where $c_{ij}$ and $q_k$ are arbitrary real numbers.

This case can be trivially handled if the original set of variables, $x$, is replaced by its log counterpart,
\[ \alpha_k= \ln{x_k} \quad \quad   k=1, \dots, n \]

Then the auxiliary vector $y$ is defined as in the previous case, avoiding duplicated entries:
\[  y_{ij}= \prod_{k=1}^{n_j} x_k^{q_k} = \prod_{k=1}^{n_j} \exp{(q_k \alpha_k)}  \]
which leads to the desired underdetermined linear system:
\[  h_i(x) = \sum_{j=1}^{m_i} c_{ij} y_{ij}  \]

The second key point is to embed the $\ln$ function in the nonlinear relationship $u=f(y)$, as follows:
\[  u_{ij}= \ln{y_{ij}}   \quad \Rightarrow \quad y_{ij}=\exp{u_{ij}} \]
which leads to the overdetermined linear system to be solved at the second step:
\[  u_{ij}=  \sum_{k=1}^{n_j} q_k \alpha_k \] % = \sum_{k=1}^{n_j} q_k \ln{x_k} \]

Once vector $\alpha$ is obtained by the factored procedure, the original unknown vector $x$ can be recovered from:
\[  {x_k}= \exp{\alpha_k} \quad \quad   k=1, \dots, n \]

\vspace{3pt}
\noindent{\bf Example 3:}

Consider the following nonlinear system in two unknowns:
\begin{eqnarray*}
p_1&=& x_1x_2 + x_1x_2^2  \\ 
p_2&=& 2x_1^2x_2 - x_1^2  
\end{eqnarray*}
which, upon introduction of four $y$ variables,
\[ y_1=  x_1x_2 \,\, ; \,\,  y_2= x_1x_2^2  \,\, ; \,\, y_3=x_1^2x_2 \,\, ; \,\, y_4=x_1^2 \]
is converted into an underdetermined linear system:
\[ 
\left(\begin{array}{c} p_1 \\ p_2 \end{array} \right)=
\underbrace{\left(\begin{array}{cccc} 1 & 1 & 0 & 0 \\ 0 & 0 & 2 & -1 \end{array} \right)}_E
\left(\begin{array}{c} y_1 \\ y_2 \\ y_3 \\ y_4 \end{array} \right)
\]

The nonlinear relationships $u=f(y)$ are,
\[  u_{i}= \ln{y_{i}} \quad\quad  i=1,2,3,4 \]
which lead to the following Jacobian:
\[ F^{-1} = \left(\begin{array}{cccc} \exp{u_1} & 0 & 0 & 0 \\ 0 & \exp{u_2} & 0 & 0 \\
0 & 0 & \exp{u_3} & 0 \\
0 & 0 & 0 & \exp{u_4} 
 \end{array} \right)
\]
and the final overdetermined system in the log variables,
\[
\underbrace{\left(\begin{array}{cc} 1 & 1 \\ 1 & 2 \\ 2 & 1 \\ 2 & 0 \end{array} \right)}_C
\left(\begin{array}{c} \alpha_1 \\ \alpha_2 \end{array} \right)=
\left(\begin{array}{c} u_1 \\ u_2 \\ u_3 \\ u_4 \end{array} \right)
\]

Once the problem is solved in the log variables, the original variables are recovered from,
\[  {x_i}= \exp{\alpha_i} \quad \quad   i=1,2,3,4 \]

For $p_1=24$ and $p_2=20$ there is a known solution at $x_1=2$ and $x_2=3$. Table \ref{t1ex3} provides the number of iterations for different starting points. Apart from taking a larger number of iterations, the NR method sometimes converges to the alternative point $x_1=31.1392$ and $x_2=0.5103$ (marked with an `A' in the table), unlike the factored scheme, which seems to converge invariably to the same point. For the farthest starting point tested (last row), the NR method does not converge in 50 iterations. On the other hand, the factored scheme is much less sensitive to the initial guess, since it converges systematically to the same point. 

\begin{table}[hbt]
\caption{Number of iterations required by the NR and the factored procedures to converge from arbitrary starting points for the nonlinear system of {Example 3}, with $p_1=24$ and $p_2=20$.}\label{t1ex3}
  \centering
  \begin{tabular}{|c|r|c|}
    \hline
    $x_0$ & NR & Factored \\ \hline
 $(1,1)$   & 7 & 6 \\
 $(1,-1)$   & 14 & 6 \\
 $(-1,1)$   & (A) 29 & 6 \\
 $(10,10)$   & 8 & 7 \\
 $(-10,-10)$   & (A) 10 & 8 \\
$(-10,10)$   & (A) 22 & 7 \\
$(-100,100)$   &  Fails & 7 \\ \hline
  \end{tabular}
\end{table}

\section{Transformation to canonical forms}\label{sec4}

There are many systems of equations which can be easily transformed into one of the canonical forms discussed above, by strategically introducing extra variables aimed at simplifying the complicating factors. Then, the nonlinear expressions defining the artificial variables are added to the original equation system, leading to an augmented one.

\vspace{3pt}
\noindent{\bf Example 4:}

Consider the following nonlinear system:
\begin{eqnarray*}
p_1&=& x_1 \sin{(x_1^2+x_2)} - x_1^2 \\ 
p_2&=& x_1^2x_2 - \sqrt{x_2}  
\end{eqnarray*}

By defining a new variable,
\[  x_3=\sin{(x_1^2+x_2)} \]
it can be rewritten in augmented canonical form, as follows:
\begin{eqnarray*}
p_1&=& x_1 x_3 - x_1^2 \\ 
p_2&=& x_1^2x_2 - \sqrt{x_2} \\
0    &=&  x_1^2+x_2 -\arcsin{x_3} 
\end{eqnarray*}

Note that in this case only the original unknowns, $x_1$ and $x_2$, need to be initialized, since $x_3$ can be set according to its definition.

\section{Extending the range of reachable solutions}\label{sec5}

When there are multiple solutions, the NR method converges to a point determined by the basin of attraction in which the initial guess is located. Some heuristics have been developed enabling the NR method to reach alternative solutions in a more or less systematic manner, such as the differential equation approach proposed in \cite{branin}, which modifies the equation structure according to the sign of the Jacobian determinant.

In this regard, the behaviour of the proposed factored algorithm is mainly affected by the computable range of $f(y)$ components. When this range does not span the whole real axis, the two-stage procedure may not be able to reach some solution points, irrespective of the starting point chosen. For instance, if $y=x^{q}$, with $q$ even, then during the solution process, the expression
\[ u=f(y)=\sqrt[q]{y}\]
will always return a positive $u$, provided $y$ is positive (or a complex with positive real component if $y$ is negative). This will be determinant of the computed $x$ values. In those cases, an easy way of allowing the factored procedure to reach alternative solution points is by considering negative ranges for $f(y)$, i.e.,
\[ u=-\sqrt[q]{y}\]

\vspace{3pt}
\noindent{\bf Example 5:}

In Example 1 it was noted that the factored procedure always converges to point B (the positive root in Fig. \ref{fex1}). 
However, if the original expression for $u_1$, 
\[  
u_1=\sqrt[4]{y_1} 
\]
is replaced by\footnote{When using Matlab, the function {\it nthroot} should be used to obtain real $q$-th roots of negative numbers, with $q$ odd. Alternatively, the use of the equivalent fractional power will provide complex roots, if needed.} 
\[  
u_1=-\sqrt[4]{y_1},  
\]
then the factored procedure always converges to point A (negative root) in fewer iterations than in Example 1. 

\noindent\underline{\hspace{6mm}}
\vspace{3pt}

A similar situation arises when the nonlinear system to be solved contains periodic functions, %such as trigonometric ones, 
since their inverses will have a limited range related to the period of the original function. For instance, if $y=\sin{x}$, then $u=\arcsin{y}$ will naturally lie in the range $-\pi/2 < u < \pi/2$. If we want to obtain $u$ values in the extended range,
\[ (q-1/2)\pi  < u < (q+1/2)\pi  \]
for any integer $q$, then we should replace $u=\arcsin{y}$ by
\[  u= q\pi + (-1)^q \arcsin{y} \]

In a similar fashion, $u=\arccos{y}$, which naturally lies in the range $0 < u < \pi$,  should be replaced by,
\[  u= (q+1/2)\pi + (-1)^q (\arccos{y}-\pi/2) \]
to obtain $u$ values in the same extended range. %,
%\[ (q-1/2)\pi  < u < (q+1/2)\pi \]

\vspace{3pt}
\noindent{\bf Example 6:}

If, in Example 2, the expression,
\[ u=f(y)= \left(  \begin{array}{c} \arcsin{y_1} \\  \arccos{y_2} \end{array} \right) 
\]
is replaced by,
\[ u= \left(  \begin{array}{c} 2\pi \\  2\pi \end{array} \right)  + \left(  \begin{array}{c} \arcsin{y_1} \\  \arccos{y_2} \end{array} \right) 
\]
then, the factored method converges, in a similar number of iterations, to the points 6.9267 and 7.2105, at a distance $2\pi$ from points A and B in figure \ref{fex2}  (in this simple example, involving a purely periodic function, this is totally an expected result, but the reader is referred to the not so trivial Example 7).

\noindent\underline{\hspace{6mm}}
\vspace{3pt}

The fact that the range of $f(y)$ adopted during the iterative process determines the solution point which is finally reached, irrespective of the starting point, is a nice feature worth investigating in the future, since by fully exploring all alternative values $u=f(y)$ one might be able to systematically reach solution points in given regions of interest, without having to test a huge number of starting points.  
The following examples illustrate this idea.

\vspace{3pt}
\noindent{\bf Example 7:}

\begin{figure}[hbt]
\begin{center}
\noindent\includegraphics[width=0.45\textwidth]{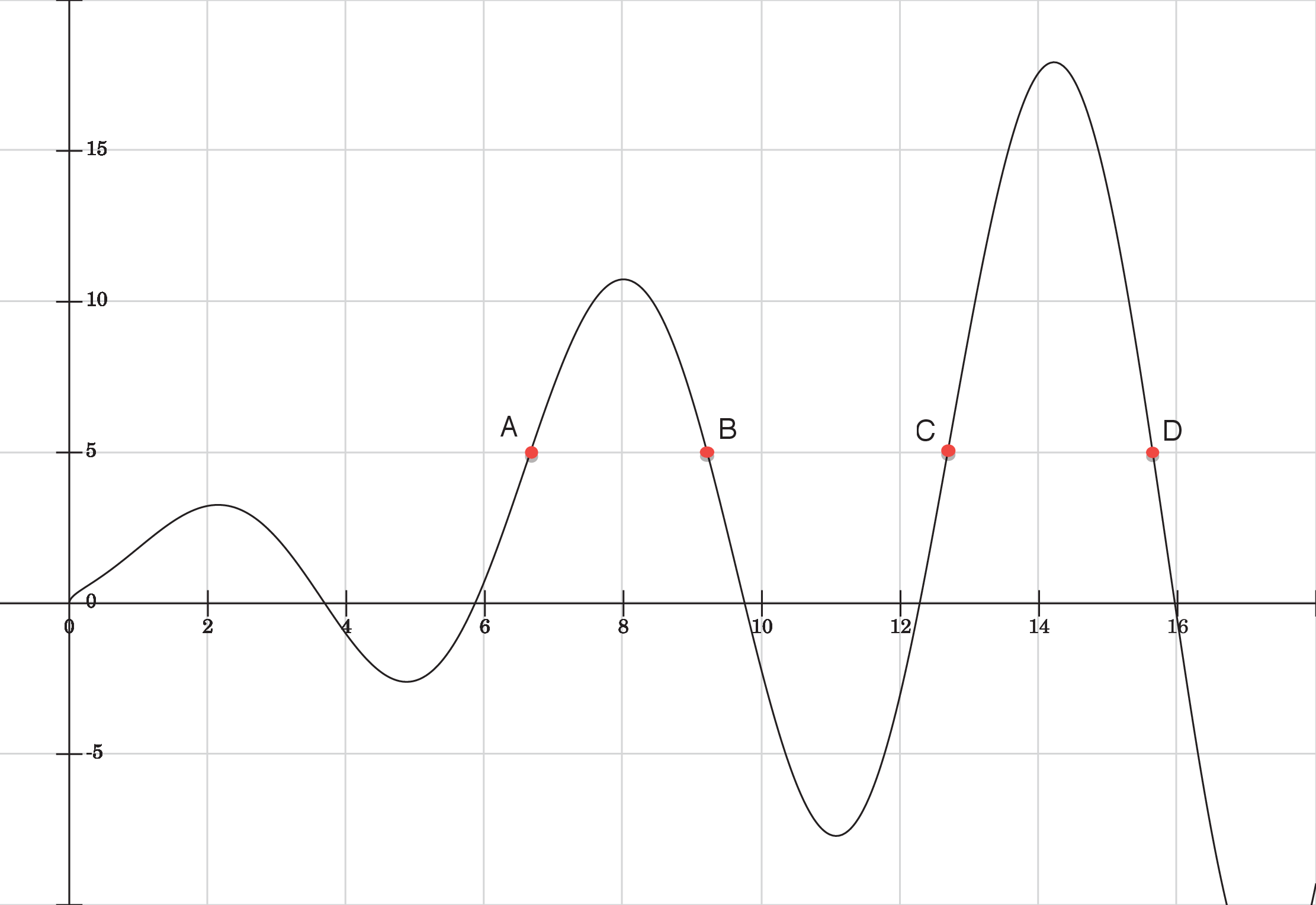}
\caption{Nonlinear function for Example 7} \label{fex7}
\end{center}
\end{figure}

Consider the extremely nonlinear non-periodic function,
\[ p= x \sin{x} +\sqrt{x} \]
which is graphically represented in figure \ref{fex7}. 

As explained in the previous examples, the following augmented system is to be factored,
\begin{eqnarray*}
p &=& x_1 x_2 +\sqrt{x_1} \\ 
0 &=& x_2 - \sin{x_1}
\end{eqnarray*}

The relevant components and relationships of the factored solution approach are,
\begin{eqnarray*}
\begin{array}{rcl}
y_1 &=& x_1 x_2 \\ 
y_2 &=& \sqrt{x_1} \\
y_3 &=& x_2 \\
y_4 &=& \sin{x_1}
\end{array}
&&
\begin{array}{rcl}
u_1 &=& \ln{y_1} \\ 
u_2 &=& \ln{y_2}  \\
u_3 &=& \ln{y_3}  \\
u_4 &=& \ln{(\arcsin{y_4})} 
\end{array}
\end{eqnarray*}

\begin{eqnarray*}
\left(\begin{array}{c} p \\ 0 \end{array}\right) &=&
\left(\begin{array}{cccc} 
1 & 1 & 0 & 0 \\
0 & 0 & 1 & -1
\end{array}\right) 
\left(\begin{array}{c} y_1 \\ y_2 \\ y_3 \\ y_4 \end{array}\right)
\end{eqnarray*}

\begin{eqnarray*}
\left(\begin{array}{cc} 
1 & 1 \\ 1/2 & 0 \\
0 & 1 \\ 1 & 0
\end{array}\right) \left(\begin{array}{c} \alpha_1 \\ \alpha_2 \end{array}\right)  &=&
\left(\begin{array}{c} u_1 \\ u_2 \\ u_3 \\ u_4 \end{array}\right) 
\end{eqnarray*}
where $\log$ variables have been introduced,
\[ \alpha_i = \ln{x_i} \quad ; \quad i=1,2   \]

The inverse of the Jacobian is in this case:
\[ F^{-1} = \left(\begin{array}{cccc} \exp{u_1} & 0 & 0 & 0 \\ 0 & \exp{u_2} & 0 & 0 \\
0 & 0 & \exp{u_3} & 0 \\
0 & 0 & 0 & \exp{u_4}\cos{(\exp{u_4})} 
 \end{array} \right)
\]

For $p=5$ no real solution exists when $x<3\pi/2$. In order to obtain real solutions for $x> 3\pi/2$ (points A, B, C and D in figure \ref{fex7}) as explained before, we have to replace, 
\[ u_4 = \ln{(\arcsin{y_4})} \]
by the more general expression,
\[ u_4 = \ln{[q\pi + (-1)^q \arcsin{y_4}]} \]
for $q>1$.

Table \ref{t1ex7} shows, for increasing values of $q$, the points to which the factored solution converges, as well as the number of iterations (starting with $x_0=q\pi$). The complex solutions provided for $q=0$ and $q=1$ (whose real component is very close to the local maximum), indicate that no real solution exists for $x< 3\pi/2$.

\begin{table}[hbt]
  \caption{Example 7: No. of iterations and solution points for $p=5$}\label{t1ex7}
  \centering
  \begin{tabular}{|c|c|c|}
    \hline
$q$  & It.   &       $x$                           \\ \hline
0  & 10  & $2.1519 + 0.5820i$     \\  
1 &  8   & $2.2158 + 1.0097i$      \\
2 &  5   &  $6.6554$ (A)                \\
3 &  5   &  $9.2097$  (B)               \\
4 &  5   &  $12.6801$ (C)               \\
5 &  4   &  $15.6411$  (D)             \\ \hline
  \end{tabular}
\end{table}

\noindent\underline{\hspace{6mm}}
\vspace{3pt}

One more example will be worked out to illustrate the capability of the factored procedure to reach different solutions, just by extending the computed range of certain $f(y)$ components. 

\vspace{3pt}
\noindent{\bf Example 8:}

Consider the 2$\times$2 system proposed by P. Boggs \cite{boggs},
\begin{eqnarray*}
-1 &=& x_1^2 -  x_2  \\ 
0 &=& x_1 - \cos{(\pi x_2/2)}
\end{eqnarray*}
which is known to have only three solutions: $(0,1)$, $(-1/\sqrt{2},3/2)$ and $(-1,2)$.

In this case, the following relationships are involved in the factored method:
\begin{eqnarray*}
\begin{array}{rcl}
y_1 &=& x_1^2 \\ 
y_2 &=& x_2 \\
y_3 &=& x_1 \\
y_4 &=& \cos{(\pi x_2/2)}
\end{array}
&&
\begin{array}{rcl}
u_1 &=& \sqrt{y_1} \\ 
u_2 &=& {y_2}  \\
u_3 &=& {y_3}  \\
u_4 &=& 2 (\arccos{y_4})/\pi 
\end{array}
\end{eqnarray*}

\begin{eqnarray*}
\left(\begin{array}{c} -1 \\ 0 \end{array}\right) &=&
\left(\begin{array}{cccc} 
1 & -1 & 0 & 0 \\
0 & 0 & 1 & -1
\end{array}\right) 
\left(\begin{array}{c} y_1 \\ y_2 \\ y_3 \\ y_4 \end{array}\right)
\end{eqnarray*}

\begin{eqnarray*}
\left(\begin{array}{cc} 
1 & 0 \\ 0 & 1 \\
1 & 0 \\ 0 & 1
\end{array}\right) \left(\begin{array}{c} x_1 \\ x_2 \end{array}\right)  &=&
\left(\begin{array}{c} u_1 \\ u_2 \\ u_3 \\ u_4 \end{array}\right) 
\end{eqnarray*}

\[ F^{-1} = \left(\begin{array}{cccc} 2{u_1} & 0 & 0 & 0 \\ 0 & 1 & 0 & 0 \\
0 & 0 & 1 & 0 \\
0 & 0 & 0 & -\pi\sin{(\pi u_4/2)/2} 
 \end{array} \right)
\]

When $u_1$ and $u_4$ are defined as above, the factored method always converges to $(0,1)$, irrespective of the starting point $x_0$. What is as important, the number of iterations is only slightly affected by $x_0$, no matter how arbitrary it is.

If we extend the computed range of $u_1$ to negative values, by replacing the original definition with,
\[ u_1 = -\sqrt{y_1} \]
then the factored method always converges to $(-1/\sqrt{2},3/2)$, irrespective of the starting point. Finally, when the ranges of both $u_1$ and $u_4$ are extended, as follows,
\[ \begin{array}{rcl}
u_1 &=& -\sqrt{y_1} \\ 
u_4 &=& 2 (2\pi-\arccos{y_4})/\pi 
\end{array}
\] 
the third solution point, $(-1,2)$, is reached for arbitrary values of $x_0$.

This kind of controlled convergence cannot be easily implemented in the NR method, characterized by rather complex attractions basins in this particular example. In fact, Newton's method may behave irregularly even when started near the solutions \cite{rheinboldt}.

Interestingly enough, if the fourth combination of ranges is selected,
\[ \begin{array}{rcl}
u_1 &=& \sqrt{y_1} \\ 
u_4 &=& 2 (2\pi-\arccos{y_4})/\pi 
\end{array}
\] 
then, the factored procedure converges to the complex point 
$(1.7174 + 0.2131i,3.9041 + 0.7320i)$
showing that no real solution exists in that region.
This issue is discussed in the next section.

\noindent\underline{\hspace{6mm}}
\vspace{3pt}

\section{Infeasibility and complex solutions}\label{sec7}

A nonlinear system $h(x)=p$ with real coefficients will have in general an undetermined number of real and complex solutions. It can be easily shown that complex solutions constitute conjugate pairs if and only if $f(y^*)=[f(y)]^*$, which happens for many functions of practical interest. When no real solutions exist we say that vector $p$ is infeasible.

The basin of attraction associated to each solution (real or complex) is characteristic of the iterative method adopted.  As discussed in section \ref{sec2}.\ref{sec2.4} and illustrated by the previous examples, the basins of attraction associated to the factored method are quite different from those of Newton's algorithm. 
Depending on which basin of attraction the initial guess $x_0$ lies in, the factored method will converge to a real or complex solution. 

A summary of the possible outcomes of the factored algorithm follows:
\begin{itemize}
\item Convergence to a real solution: for feasible cases, starting from a real $x_0$ close enough to a solution, the factored scheme eventually converges to the nearby real solution. This is favoured by the introduction of the first step, aimed at minimizing the distance between successive intermediate points. Note however that, even if both $x_0$ and the solution are real, depending on the domain of existence of the nonlinear elementary functions $f({y})$, %and $\tilde{F}^{-1}$, 
the two-stage procedure may give rise to complex intermediate values $x_k$ during the iterative process. A real solution  (or, more precisely, a solution with negligible imaginary component) can also be reached starting from a complex $x_0$.
\item Convergence to a complex solution: for infeasible cases the factored method can only converge to a complex solution, provided the right $x_0$ is chosen. Notice that complex solutions cannot be reached if $x_0$ is real and the domains of $f({y})$ and ${F}^{-1}$ span the entire real axis. In those cases, selecting a complex initial guess is usually helpful to allow the algorithm to reach a complex solution. For feasible cases, if $x_0$ is real but far from a real solution, the factored iterative procedure might converge to a complex point. 
\item Inability to converge: Like any iterative algorithm, the factored scheme may numerically break down owing to unbounded values provided by the elemental functions $f(\cdot)$ or their inverse $f^{-1} (\cdot)$, which is the case for instance when an asymptote is approached. It can also remains oscillating, either periodically or chaotically, which typically happens when $x_0$ is real, no real solution exists and complex intermediate values cannot pop up. 
\end{itemize}

The following examples illustrate the behaviour of the factored scheme as well as the nature of the complex solutions obtained when solving infeasible cases.

\vspace{3pt}
\noindent{\bf Example 10:}

Let us reconsider Example 2 with $p=1.5$, which constitutes an infeasible case for which the NR fails to converge to a meaningful point. The two-stage procedure always converges in a reduced number of iterations to the same complex value  
(or its conjugate), as shown in table \ref{t1ex10} for the same real starting points as in Example 2 (the same happens with other real starting points). 

\begin{table}[hbt]
  \caption{Example 10: No. of iterations and solution points for $p=1.5$}\label{t1ex10}
  \centering
  \begin{tabular}{|r||c|c|}
    \hline
    $x_0$ & Iter & $x$ \\ \hline
10      & 8 & $0.7854 - 0.3466i$ \\
5      & 5 & $0.7854 + 0.3466i$  \\
1      & 8 & $0.7854 - 0.3466i$   \\
0      & 5 & $0.7854 + 0.3466i$  \\
$-1$ &  5 & $0.7854 + 0.3466i$  \\
$-5$  & 6  & $0.7854 - 0.3466i$   \\
$-10$ & 5 & $0.7854 - 0.3466i$ \\
    \hline
  \end{tabular}
\end{table}

However, if we prevent $f(y)$ from taking complex values by restricting $y$ to its real domain,  
then the two-stage procedure remains oscillating around real values, much like the NR method.

We can use this example to see what happens if $p$ is further increased. Table \ref{t2ex10} presents the number of iterations and the solution points for increasing values of $p$ (starting from $x_0=0$). The maximum value for which there is a feasible (real) solution is $p=1.4142$ (critical  point). For larger $p$ values the factored approach converges to complex values with the same real component and increasing absolute values. Eventually, for $p=4.204$ the two-step procedure breaks down, which means that $x_0=0$ is not in the basin of attraction for this infeasible value of $p$. 

\begin{table}[hbt]
  \caption{Example 10: No. of iterations and solution points for feasible and infeasible values of $p$, starting from $x_0=0$ 
  }\label{t2ex10}
  \centering
  \begin{tabular}{|r||c|c|}
    \hline
    $p$ & Iter & $x$ \\ \hline
1.4        & 7 & $0.6435$ \\
1.4142   & 12 & $0.7810$  \\ \hline
1.4143   & 10 & $0.7854 + 0.0111i$   \\
1.5        & 5 & $0.7854 + 0.3466i$  \\
%2           & 4 & $0.7854 + 0.8814i$ \\
2.5        &  5 & $0.7854 + 1.1711i$  \\
%2.768    & 5  & $0.7854 + 1.2919i$   \\ \hline
%2.769    & 5 & $0.7854 + 1.2923i$ \\
3           & 5 & $0.7854 + 1.3843i$ \\
4.203    & 10 & $0.7854 + 1.7528i$ \\ \hline
4.204    & \multicolumn{2}{c|}{Fails} \\
    \hline
  \end{tabular}
\end{table}

Experimental results show that the complex values to which the factored method converges in infeasible cases are not  arbitrary. In this example, evaluating the nonlinear function at the real component of the solution point (0.7854) yields $p_r=1.4142$ (point C in figure \ref{fex2}). This lies almost exactly on the maximum of the function ($p_M=\sqrt{2}$).  

In Example 1, if we set $p=-0.2$, which is infeasible in the real domain, the factored method converges from nearly arbitrary real starting points to the complex value $0.8090 + 0.2629i$. Taking the real part of the solution yields $p_r=-0.1011$
(point C in figure \ref{fex1}), which is indeed very close to the minimum of the function $x^4-x^3$ ($p_m=-0.1055$ at $x=0.75$).

\noindent\underline{\hspace{6mm}}
%\vspace{3pt}

\vspace{3pt}
\noindent{\bf Example 11:}

\begin{figure}[hbt]
\begin{center}
\noindent\includegraphics[width=0.45\textwidth]{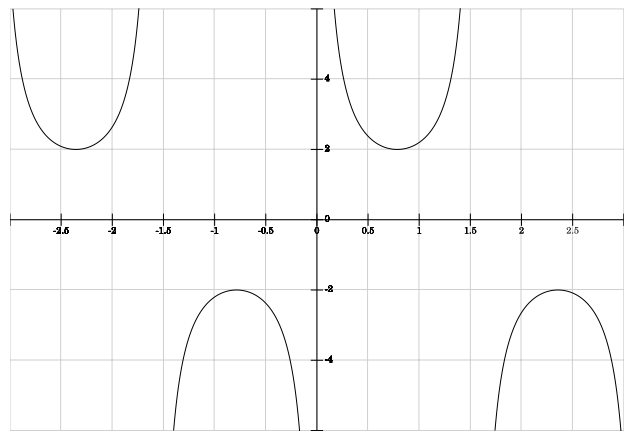}
\caption{Periodic function for Example 10} \label{fex11}
\end{center}
\end{figure}

Let us consider the period function,
\[ p= \tan{x} - \tan{(x-\pi/2)} \]
which is infeasible for $-2<p<2$ and has an infinite number of solutions otherwise (see figure \ref{fex11}). 

This example is very similar to Example 2, except for the way the auxiliary variables are defined:
\[ y_1 = \tan{x} \quad ; \quad y_2= \tan{(x-\pi/2)} \]

In this case, both the elementary nonlinear functions $f(y)$, 
\[ u_1 = \arctan{y_1} \quad ; \quad u_2= \pi/2+\arctan{y_2}
\]
and the Jacobian functions,
\[ F^{-1} = \left(\begin{array}{cc} 1+\tan^2{u_1} & 0 \\ 
                                                                                   0 & 1+\tan^2{(u_2-\pi/2)} 
 \end{array} \right)
\]
are well defined all over the real axis (except for the asymptotes of the $\tan$ functions).

Therefore, starting from a real value $x_0$, the factored procedure (and also the NR scheme) will always remain in the real domain, unlike in Example 2 containing the functions $\arcsin$ and $\arccos$.

With $p=3$, starting from $x_0=1$, the two-stage algorithm converges in 5 iterations to 1.2059, whereas the solution point 0.3649 is reached, also in 5 iterations, when $x_0=-1$.

For the infeasible value $p=1.9$, the algorithm remains oscillating for any real starting point. However, with the complex starting value $x_0=1+i$, it converges to the complex solutions shown in table \ref{t1ex11} for different infeasible values of $p$. 

\begin{table}[hbt]
  \caption{Example 11: No. of iterations and complex solution points starting from $x_0=1+i$}\label{t1ex11}
  \centering
  \begin{tabular}{|c||c|c|}
    \hline
    $p$ & Iter. & $x$ \\ \hline
$1.9$  & 6 & $0.7854 + 0.1615i$ \\
$1.5$   & 4 & $0.7854 + 0.3977i$  \\
$1$  &   4  & $0.7854 + 0.6585i$  \\ \hline
  \end{tabular}
\end{table} 

Note that the real component of the solution is always the same (0.7854). This value is exactly $\pi/4$, a local minimum of the nonlinear function.

\noindent\underline{\hspace{6mm}}
\vspace{3pt}

The above examples suggest that, when the two-stage method fails to find a real solution, it tries to reach a complex point whose real component is in the neighborhood of an extreme point of the function, which is indeed a nice feature. This is surely a consequence of the strategy adopted in step 1, which attempts to minimize the distance from new solution points to the previous ones, but this conjecture deserves a closer examination which is out of the scope of this preliminary work. 

In any case, it is advisable to adopt complex starting points, particularly when infeasibility is suspected, in order to facilitate complex solutions being smoothly reached.

\section{Critical points}\label{sec8}

Feasible (real) and infeasible (complex) solution points can be obtained as long as the Jacobian matrix ($H_k=EF_k^{-1}C$) remains nonsingular throughout the iterative process. This applies to both the NR and the factored methods. Even though critical points are theoretically associated with singular Jacobians, in practice they can also be safely approached so long as the condition number of the Jacobian is acceptable for the available computing precision and convergence threshold adopted. However, the convergence rate tends to significantly slow down near critical points.

It is worth noting that, in the surroundings of a critical point, it may be preferable to deal with the augmented equation (\ref{FOOCauxnoyc}) rather than solving the more compact one (\ref{newx0}), which assumes that $\mu=0$. 

\vspace{3pt}
\noindent{\bf Example 12:}

Let us consider again the periodic nonlinear function of Example 11:
\[ p= \tan{x} - \tan{(x-\pi/2)} \]
which has an infinite number of critical points for $p=2$. Table \ref{t1ex12} provides the number of iterations and the solution points reached by both the NR and factored methods for different starting points (like in previous examples the convergence threshold is $||\Delta x||_1< 0.00001$). As can be seen, the factored method is much more robust against the starting point choice, and converges always to the same critical point ($\pi/4$).

\begin{table*}[hbt]
  \caption{Example 12: No. of iterations and solution points for $p=2$ and $p=2.1$}\label{t1ex12}
  \centering
\begin{tabular}{|r||c|r||c|r||c|r||c|r|}
  \hline
             & \multicolumn{4}{c||}{$p=2$} & \multicolumn{4}{c|}{$p=2.1$} \\ \cline{2-9}
  $x_0$ & \multicolumn{2}{c||}{Factored} & \multicolumn{2}{c||}{NR} & \multicolumn{2}{c||}{Factored} & \multicolumn{2}{c|}{NR} \\ \cline{2-9}
     &  Iter & $x$ &  Iter & $x$ & Iter & $x$ &  Iter & $x$ \\ \hline
5 & 16  & 0.7854 & 23 & $-178.2854$ & 5  & 0.6305 & 8 & $-37.0686$  \\
3 & 15  & 0.7854 & 25 &  101.3164 & 6  & 0.6305 & 13 & 4.0819  \\
1.5 & 16  & 0.7854 & 19 &  0.7854  & 6  & 0.9403 & 9 &  0.9403  \\
$-1.5$ & 16  & 0.7854 & 20 & $-2.3562$  & 6  & 0.6305 & 12 & $-2.2013$  \\
$-3$ & 15  & 0.7854 & 18 & $-2.3562$ & 6  & 0.9403 & 8 & $-2.5111$   \\
$-5$ & 16  & 0.7854 & 17 &  $-5.4978$  &  5 & 0.9403 & 6 & $-5.3429$   \\ \hline
\end{tabular}
\end{table*}

In the neighbourhood of the critical point the number of iterations is significantly reduced, but the factored procedure still performs much better than the NR method. For the sake of comparison, Table \ref{t1ex12} also shows the values corresponding to $p=2.1$.

\noindent\underline{\hspace{6mm}}
\vspace{3pt}

In addition, numerical problems may arise during the factored solution procedure if, by a numerical coincidence, any element of the diagonal matrix $F^{-1}$ becomes null or undefined. This risk can be easily circumvented by preventing diagonal elements of $F^{-1}$ from being smaller than a certain threshold or abnormally large.

\section{Large-scale test cases} \label{lstc}

This section provides test results corresponding to the nonlinear systems arising in the steady-state analysis of electrical power systems, known in the specialized literature as the ``power flow'' problem \cite{stott}. The results presented herein are representative of many other engineering problems characterized by very sparse equations with network or graph-based structure, in which state variables associated to nodes are linked directly only to their neighbours in the graph \cite{sandberg}.

In the power flow problem, the active and reactive power injections specified at each bus $i$, $P_i^{sp}$ and $Q_i^{sp}$ respectively, are linearly related with unknown powers %($P_{ij}$ and $Q_{ij}$) 
flowing through branches linking $i$ with its neighbours $j\in i$, as follows:
\begin{eqnarray}\label{ntb}
P_i^{sp} &=& \sum_{j\in i} P_{ij}  \\
Q_i^{sp} &=& \sum_{j\in i} Q_{ij} \nonumber
\end{eqnarray}
In turn, branch power flows are nonlinearly related to nodal state variables (voltage magnitudes, $V_i$, and phase angles $\theta_i$):
\begin{eqnarray}\label{btx}
P_{ij} &=& g_{ij}V^2_i-g_{ij}V_iV_j\cos{\theta_{ij}}-b_{ij}V_iV_j\sin{\theta_{ij}}  \\
Q_{ij} &=& -(b_{sh}+b_{ij})V^2_i+b_{ij}V_iV_j\cos{\theta_{ij}}-g_{ij}V_iV_j\sin{\theta_{ij}} \nonumber
\end{eqnarray}
where $\theta_{ij}=\theta_i-\theta_j$ and $g_{ij}$, $b_{ij}$ and $b_{sh}$ represent branch parameters \cite{stott}. For a network with $N$ buses the phase angle of a given node is arbitrarily taken as a reference ($\theta_s=0$), while the active power injection at this same node (termed slack bus) is not specified to account for unknown network losses. Replacing the pair (\ref{btx}) into (\ref{ntb}) leads therefore to a set of $n=2N-1$ nonlinear equations in the same number of unknowns, which are customarily solved by means of the NR method. In practice, the voltage magnitude rather than the reactive power is specified at generator buses, accordingly reducing the size of the equation system. Using the so-called flat-start profile ($V_i=1$, $\theta_i=0$), the nonlinear power flow equations are solved typically in 3 to 5 iterations by Newton's method, depending on the network loading level and ill-conditioning. 

As explained in \cite{flf}, the natural choice for the intermediate vector $y$ in this particular application is
\begin{equation}\label{eqn6}
y=\{ U_i, K_{ij}, L_{ij} \}
\end{equation}
where $U_i = V_i^2$ for every node $i$ and, for each branch connecting nodes $i$ and $j$, the following pair of variables is adopted:
\begin{eqnarray}\label{eqn1}
   K_{ij} &=& V_iV_j\cos \theta_{ij} \\
   L_{ij} &=& V_iV_j\sin \theta_{ij}
\end{eqnarray}
Vector $y$ then comprises $m=2b+N$ variables ($b$ being the number of branches). It can be easily checked that, when expressed in terms of $y$, the power flow equations become linear (the reader is referred to \cite{flf} for the details).

Table \ref{tlfsol} compares the number of iterations and solution times (miliseconds) for both the proposed factored procedure and Newton's method, when applied to several benchmark systems ranging from 30 to over 3000 nodes (i.e. some 6000 nonlinear equations) available in \cite{matpower}. The convergence threshold adopted is $||\Delta p_k||_\infty<0.001$.

\begin{table}[hbt]
  \caption{No. of iterations and solution times (ms.) for the factored and NR methods.} \label{tlfsol}
  \centering
\begin{tabular}{|r||c|c||c|c|} \hline
       & \multicolumn{2}{c||}{Factored} & \multicolumn{2}{c|}{NR} \\ \cline{2-5}
$N$ &  Iter. & Time &  Iter. & Time \\ \hline
30     &  2	& 1.36	& 3    & 1.97  \\
39     &  3	& 1.87	& 4    & 2.45  \\
57     &  3	& 2.39	& 3    & 2.56  \\
%118   &  3	& 3.7	& 3    & 3.6    \\
300   &  3	& 8.4	& 4    & 10.01 \\
2383	&  3  & 64.21	& 4    & 73.54 \\
%2736	&  4  & 89.46	& 5    & 111.22 \\
2737	&  4  & 102.9	& Fails &  $-$ \\
%2746	&  5  & 112.26	& 6    &  129.94 \\
3012	&  4  & 99.15	& Fails & $-$ \\
3120 &  4	& 102.02	& 5    & 126.16 \\
\hline
\end{tabular} \\
\end{table}

Solution times have been obtained on an Intel Core i7 processor (2.8 GHz, 8 GB of RAM) by running a prototype Matlab implementation. Every effort has been made to optimize the code in both cases, by exploiting available sparsity features and avoiding duplicated operations (e.g., the Jacobian $H_k$ in Newton's method is obtained almost free in this particular application, as a byproduct of the computations involved in the mismatch vector $\Delta p_k$).

As can be seen, the factored method always converges and takes one iteration less than NR in most cases. In scenarios representing highly loaded conditions (i.e. more nonlinearity), where the customary flat start is not so close to the solution point, the factored algorithm tends to save two or more iterations with respect to NR method \cite{flf}.

\begin{figure}[hbt]
\begin{center}
\noindent\includegraphics[width=0.45\textwidth]{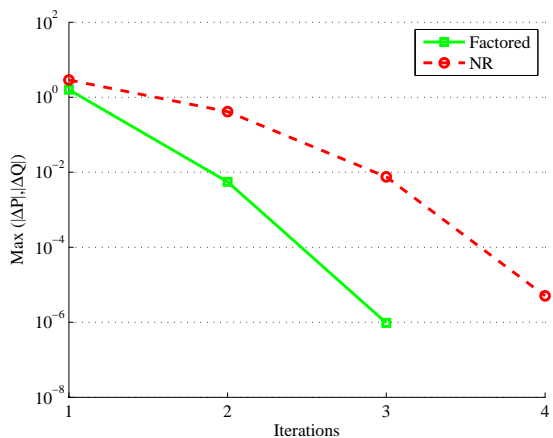}
\caption{Convergence rates of NR and factored methods for the 300-bus system} \label{conv300}
\end{center}
\end{figure}

Further insight into the local convergence rate of the factored method can be gained from Fig. \ref{conv300}, where the evolution of $||\Delta p_k||_\infty$ when solving the 300-bus system from flat start is represented, both for Newton's and factored methods.

\section{Conclusion}

In this article a new iterative procedure to reliably solve nonlinear systems of equations is presented and illustrated with several examples. The basic idea is to unfold the original system with the help of auxiliary variables so that elementary nonlinearities with explicit inverse can be individually handled. At each iteration the algorithm first solves a trivial least-distance problem; then a Newton-like step is computed.
For realistic problems, the overall computational effort required by the two steps should be of the same order or even lower than that of the standard NR scheme, provided suitable sparsity techniques are adopted.

Practical experience so far, including large sets of nonlinear equations arising in power systems problems, shows that the proposed method clearly outperforms the NR method near the solution point, while its behaviour, far away from NR typical basins of attraction, is less erratic and more controllable. Future efforts should be directed to more thoroughly investigating its behaviour in terms of global convergence and capability to deal with infeasible cases.

Applying the factored scheme to the development of globally-convergent and/or higher-order iterative methods, as well as other related problems such as nonlinear programming and the solution of nonlinear ordinary differential equations, also seems to be a promising venue of research. 

\section*{Acknowledgment}

The author acknowledges Walter Vargas for providing the numerical results in Table \ref{tlfsol} and Catalina G\'omez-Quiles for Fig. \ref{conv300}.

%\vfill
%\newpage

\vfill

\end{document}